\documentclass[preprint]{imsart}
\RequirePackage[OT1]{fontenc}
\RequirePackage{amsthm,amsmath}
\RequirePackage{natbib}
\usepackage{color}
\RequirePackage[colorlinks,citecolor=blue,urlcolor=blue]{hyperref}
\RequirePackage{soul}
\usepackage{amsfonts}
\usepackage{dcolumn}
\usepackage{textcomp}
\usepackage{longtable}
\usepackage{verbatim}
\usepackage{latexsym}
\usepackage{fancyvrb}
\usepackage{multirow}
\usepackage[linesnumbered,ruled,vlined]{algorithm2e}
\usepackage{url}
\usepackage{verbatim}
\usepackage{algorithmic}
\usepackage{comment}
\usepackage{natbib}
\usepackage{float}
\usepackage{graphicx}

\arxiv{math.PR/0000000}
\startlocaldefs
\numberwithin{equation}{section}
\theoremstyle{plain}

\endlocaldefs
\newtheorem{definition}{Definition}

\newtheorem{theorem}{Theorem}

\begin{document}
\begin{frontmatter}

\title{Vanilla Lasso for sparse classification under single index models}
\runtitle{Vanilla Lasso for sparse classification under generalized single index models}

\begin{aug}
\author{\fnms{Jiyi} \snm{Liu }\ead[label=e1]{nbljy111@163.com}}
\and
\author{\fnms{Jinzhu} \snm{ Jia}\ead[label=e2]{jzjia@math.pku.edu.cn}}
%

\runauthor{J. Liu and J. Jia.}
\affiliation{School of Mathematical Sciences and Center of Statistical Science,\\ Peking University} 
\address{School of Mathematical Sciences and Center of Statistical Science, LMAM, LMEQF\\
         Peking University\\
         Beijing 100871, China \\
\printead{e1}
\phantom{E-mail:\ } \printead*{e2}\\}
%
%
%
\end{aug}
\maketitle

\begin{abstract}
This paper study sparse classification problems.  We show that under single-index models, vanilla Lasso could give good estimate of unknown parameters. With this result, we see that even if the model is not linear, and even if the response is not continuous, we could still use vanilla Lasso to train classifiers. 
Simulations  confirm that vanilla Lasso could be used to get a good estimation when data are generated from a logistic regression model.
\end{abstract}
\end{frontmatter}

\section{Introduction}

Classification problem is  important  in a lot of fields such as pattern recognition, bioinformatics etc. There are a few classic classification methods, including LDA (linear discriminant analysis), logistic regression, naive bayes,  SVM (support vector machine). When the number of features (or predictors) denoted by $p$ is fixed, under some regularity conditions, LDA is proved to be optimal (see standard statistical text books such as \cite{anderson1958introduction}). Logistic regression is very popular because it has no distribution assumption on $X$. It has been proved that LDA is equivalent to least squares \cite{fisher1936use}. Using this connection, in fact, one could solve the LDA problem by directly using the vanilla $\ell_1$ penalized least squares, i.e. the Lasso (\cite{tibshirani1996regression}). Using Lasso to solve sparse LDA has already been proposed by  \cite{mai2012direct}, in which they showed that under irrepresentable condition and some other regularity conditions, the Lasso could select the important features (or predictors) for linear discriminant analysis.  $\ell_1$ penalized logistic regression has been widely used for high-dimensional classification problems \citep{koh2007interior,ravikumar2010high}. 

LDA and logistic regression model are parametric models. Both models require a few assumptions on the data collected. To make less assumptions, non-parametric or semi-parametric models could be used.  To control model complexity, single index model is a great choice for semi-parametric models \citep{powell1989semiparametric,klein1993efficient,ichimura1993semiparametric}. A single-index model is defined as follows.
\[E(y|x) = F(x^T\beta),\]
where $F$ is an unknown function. $F$ is usually estimated via Nadaraya-Watson nonparametric estimator
\[
\tilde F(x_i^T\beta) = \sum_{j\neq i} y_jK\left(\frac{x_i^T\beta - x_j^T\beta}{h}\right) /  \sum_{j\neq i} K\left(\frac{x_i^T\beta - x_j^T\beta}{h}\right).
\]
Once the form of $F$ is given, a maximum (quasi)likelihood estimation for $\beta$ could be given via maximizing the following quasi-likelihood:
\[\sum_{i=1}^n(y_i \log(\tilde F(x_i^T\beta)) + (1-y_i)\log(1-\tilde F(x_i^T\beta))).\]
Nonlinear least squares or Nonlinear weighted least square could also be used to estimated $\beta$.
Nonlinear weighted least squares estimator is defined as follows.
\[
\hat\beta = \arg\min_{\beta}\frac{1}{n}\sum_{i=1}^n w(x_i) [y_i - \tilde F(x_i^T\beta)]^2,
\]
when $w(x_i) = 1$, for all $i=1,2,\ldots,n$, it is a nonlinear least squares estimator.

For high-dimensional classification problems, a very demand is to select very few important features.  AIC or BIC could be used. Their convex relaxation -- $\ell_1$ constrained optimization could also be used. 
In this paper, we propose a new feature selection method for high-dimensional classification problems. We apply vanilla Lasso to this high-dimensional classification problems. We show that if data are from a  single index model \citep{hardle1993optimal,cui2009generalized} with $F$ unknown, then vanilla Lasso could be used to select important features. 

The rest of the paper is organized as follows. In Section \ref{sec:MR}, we give the main results of the paper. Section \ref{sec:proof} provides proofs. We use a simulate study to support our theories in Section \ref{sec:simu} and finally we conclude in Section \ref{sec:conclusion}.

\section{Vanilla Lasso for single index models}
\label{sec:MR}
Suppose  observations $(x_i^T,y_i)\in \mathbb R^{p} \times \mathbb \{-1,1\}, i=1,2,\ldots,n$ i.i.d.\ follow some unknown distribution.  $y_i$ and $x_i$ have the following single index model:
\begin{equation}
\label{eqn:single-index}
E(y_i | x_i) = F(x_i^T\beta),
\end{equation}
where $F(\cdot)$ is an unknown function and $\beta$ is the unknown parameter that could be used to construct a classifier. If $F(x) = \frac{e^x}{1+e^x}$, then we have logistic regression model; if $F(x) = P(Z\leq x)$, where $Z$ is a standard Gaussian random variable, then we have a probit model. Our goal is to estimate the unknown coefficient $\beta$ and then construct a classifier. For both logistic regression and probit model the classifier is $sign(x_i^T\beta)$. \cite{plan2013robust} first considered this problem as a one-bit compressed sensing problem. They proposed to estimate $\beta$ by maximizing $\sum_{i=1}^n y_i x_i^T \beta$ and gave a few constraints on the solution of $\beta$.  We argue that a more natural objective function like minimizing least squares could be used. We state our main results as follows.

\begin{theorem}
\label{thm:main}
Let ${\beta^*}\in \mathbb R^n$ with $\|\beta^*\|_1\le \sqrt{s},\|\beta^*\|_2=1$. The single index model \eqref{eqn:single-index} is satisfied. Let $\lambda = {E}[F({Z}){Z}]$, where ${Z}$ is a standard  Gaussian random variable. Denote $\hat \beta$ as the solution of following Lasso problem:
\[\hat \beta  = \arg\min_\beta \sum_{i=1}^n(y_i - x_i^T\beta)^2, \ \ \ s.t.\   \|\beta\|_1 \leq s.\]
Suppose $n\ge C\delta^{-2}s\log(\frac{\lambda^2p}{2s})$, then with probability $P\ge 1-13e^{-c\delta \lambda\sqrt{n}}$:
\[\left\|\frac{\hat \beta}{||\hat \beta||_2} - \beta^*\right\|_2 \le \delta .\]
\end{theorem}
\textbf{Remark 1.} \cite{plan2015generalized} has a similar result. Here we use a quite different proof skill to prove this result. Our proof depends on standard concentration inequalities.

\textbf{Remark 2.} It is easy to see that parameter $\beta$ is not identifiable, because $F(\cdot)$ is an unknown function, any scaling factor could be absorbed into $F(\cdot)$. It is reasonable that Theorem \ref{thm:main} controls the error for scaled $\beta.$ What is worthy noticing is that the restrict of original signal ${\beta^*}$ is basically $\frac{||{\beta^*}||_1}{||{\beta^*}||_2}\le \sqrt{s}$, which is called effective sparsity in previous work. We can replace ${\beta^*}$ with $\frac{{\beta^*}}{\|{\beta^*}\|_2}$ to make it normalized.
\section{Proofs}
\label{sec:proof}
We denote $\lambda = {E}[F(Z)Z]$, where ${Z}$ is a standard  Gaussian random variable. In any where in the following context, we have to ensure that $\lambda >0$. Without loss of generality, we assume $\lambda<1$ if not mentioned.
\newtheorem{theorum}{Theorum}[section]


We can combine the following statements to prove our main result stated in Theorem \ref{thm:main}.
\begin{theorem}
\label{thm:1.2}
Let $\beta^* \in \mathbb R^n$ with $ ||\beta^*||_1\le \sqrt{s},||\beta^*||_2=1$, and the single index model \eqref{eqn:single-index} is satisfied. Let $\hat\beta\left(\frac{2}{\lambda},1\right)$  to be the solution of the following constrained optimization problem. 
\[\hat\beta\left(\frac{2}{\lambda},1\right)  = \arg\min_\beta \|{y}-{X\beta}\|_2\hspace{1cm}s.t.\ \|\beta\|_1\le \frac{2\sqrt{s}}{\lambda} ,\|\beta\|_2= 1.\]
Suppose $n\ge C\delta^{-2}s\log(\frac{\lambda^2p}{2s})$, 
then with probability $P\ge 1-5e^{-c\delta \lambda\sqrt{n}}$, we have
\[\left\|\hat\beta\left(\frac{2}{\lambda},1\right) - \beta^*\right\|_2 \le \delta \]
\end{theorem}

\begin{theorem}
\label{thm:1.3}
Let $\beta^* \in \mathbb R^n$ with $ \|\beta^*\|_1\le \sqrt{s},\|\beta^*\|_2=1$, and the single index model \eqref{eqn:single-index} is satisfied. Let 
Let $\hat\beta\left(\frac{2}{\lambda},1\right)$ and  $\hat\beta\left(1,k\right)$ to be the solution of the following constrained optimization problems for k $\in [\lambda-\varepsilon,\lambda+\varepsilon]$ and $\varepsilon \le \frac{\lambda}{2}$. 
\[\hat\beta\left(\frac{2}{\lambda},1\right)  = \arg\min_\beta ||{y}-{X\beta}||_2\hspace{1cm}s.t.\ \|\beta\|_1\le \frac{2\sqrt{s}}{\lambda} ,\|\beta\|_2= 1.\]

\[\hat\beta\left(1,k\right)  = \arg\min_\beta ||{y}-{X\beta}||_2\hspace{1cm}s.t.\ \|\beta\|_1\le \sqrt{s},\|\beta\|_2= k.\]

Then when $n>C\delta^{-2}\lambda^{-2}s\log(\frac{2(\lambda-\varepsilon)^2p}{s})$, with probability $P\ge 1-6e^{-c\delta \lambda\sqrt{n}}$, 
\[\left\|\hat\beta\left(\frac{2}{\lambda},1\right)-\frac{\hat\beta(1,k)}{k}\right\|_2 \le \delta \]
\end{theorem}

\begin{theorem}
\label{thm:1.4}
Let $\beta^* \in \mathbb R^n$ with $ \|\beta^*\|_1\le \sqrt{s},\|\beta^*\|_2=1$, and the single index model \eqref{eqn:single-index} is satisfied. Let 
Let $\hat \beta $ to be the solution of the following constrained optimization problem. 
\[\hat \beta  = \arg\min_\beta \sum_{i=1}^n(y_i - x_i^T\beta)^2, \ \ \ s.t.\   \|\beta\|_1 \leq s.\]
Then with probability $P\ge 1-e^{-\frac{\varepsilon^2n}{32}}-e^{-\frac{(1-\lambda-\varepsilon)^2\varepsilon^2}{8}n}$,
\[\lambda-\varepsilon \le||\hat{\beta}||_2 \le \lambda +\varepsilon \]
\end{theorem}
Here, we show how these three combined can prove Thm 1. We abbreviate $\hat{\beta_1}$ for $\hat\beta\left(\frac{2}{\lambda},1\right)$ and set $\varepsilon=\min\{\frac{\lambda}{2},\frac{1-\lambda}{2}\}.$
 From theorem \ref{thm:1.4}, with $P\ge 1-e^{-\frac{\varepsilon^2n}{32}}-e^{-\frac{(1-\lambda-\varepsilon)^2\varepsilon^2}{8}n}$, the solution $\hat{\beta_0}$ in \eqref{eqn:single-index} satisfies $\lambda-\varepsilon \le\|\hat{{\beta_0}}||_2 \le \lambda +\varepsilon$. In this case, the condition in the theorem \ref{thm:1.3} holds. Combined with the result of theorem \ref{thm:1.2} setting $k=\frac{1}{\lambda-\varepsilon}$, we conclude
\begin{equation}
\label{finalstep}
\left\|{\beta}-\frac{\hat {\beta_0}}{\|\hat {\beta_0}\|_2}\right\|_2\le \|{\beta}-\hat{{\beta_1}}\|_2+\|\hat{{\beta_1}}-\frac{\hat {\beta_0}}{\|\hat {\beta_0}\|_2}\|_2\le 2\delta
\end{equation}
The probability for this assertion is
\begin{equation*}
\begin{aligned}
P&\ge 1-(1-P_1)-(1-P_2)-(1-P_3)\\
&\ge 1-5e^{-c\delta \lambda\sqrt{n}}-6e^{-c\delta \lambda\sqrt{n}}-e^{-\frac{\varepsilon^2n}{32}}-e^{-\frac{(1-\lambda-\varepsilon)^2\varepsilon^2}{8}n}\\
&= 1-11e^{-c\delta \lambda\sqrt{n}}-e^{-\frac{\varepsilon^2n}{32}}-e^{-\frac{(1-\lambda-\varepsilon)^2\varepsilon^2}{8}n}
\end{aligned}
\end{equation*}
Notice in \eqref{ep} we set $\varepsilon=\min\{\frac{\lambda}{2},\frac{1-\lambda}{2}\}$
\[P\ge 1-13e^{-c\delta \lambda\sqrt{n}}\]
The lower bound of $n$ in Theorem \ref{thm:1.2} and \ref{thm:1.3} can both be rewritten as $C\delta^{-2}s\log(\frac{\lambda^2p}{2s})$ under this circumstances, and thus, this bound suffices to conclude our results. Finally, replace $\delta$ with $\frac{\delta}{2}$, result in \eqref{finalstep} is exactly what we want, and thus our proof is completed.
We now prove Theorems  \ref{thm:1.2}, \ref{thm:1.3} and \ref{thm:1.4}.
\subsection{Proof of Thm \ref{thm:1.2}}

\label{secpro2}
Here we give the proof of Theorem \ref{thm:1.2}.\\

Denote $f({\beta'})=\frac{1}{n}(\sum\limits_{i=1}^n y_i\langle {x_i},{\beta'}\rangle -\frac{1}{2}\sum\limits_{i=1}^n \langle {x_i},{\beta'}\rangle^2)$.
According to the main theorem in \cite{plan2013robust}, set K=\{$\beta\in \mathbb R^n: \|\beta\|_2=1,\|\beta\|_1\le j\sqrt{s}$\} there, and we can conclude that when $\|\beta-\beta'\|_2\ge \delta$, with $P\ge 1-4e^{-2\gamma^2}$, 
\begin{equation}
\label{plan2013}
\frac{1}{n}(\sum\limits_{i=1}^n y_i\langle {x_i},{\beta}\rangle-\sum\limits_{i=1}^n y_i\langle {x_i},{\beta'}\rangle) \ge \frac{\lambda\delta}{2}-\frac{\sqrt{Cj^2s\log(\frac{2p}{j^2s})}+4\gamma}{\sqrt{n}}
\end{equation}
Noticing that $\sum\limits_{i=1}^n \langle {x_i},{\beta'}\rangle^2 \sim ||\beta'||_2^2 \chi^2(n)$, which belongs to sub-exponential distribution defined as follows.  
\begin{definition}[Sub-exponential]
A random variable $X$ is sub-exponential with  parameter $(\sigma,b)$, if  for all $|\lambda| < 1/b$,
\begin{equation}\label{eqn:sub-exp}
E\left[\exp\left[\lambda\left(X - E(X)\right)\right]\right] \leq \exp\left(\frac{1}{2}\sigma^2\lambda^2\right)
\end{equation}
\end{definition}
It is not hard to check that the parameter $(a,b)$ for $\chi^2(n)$ can be chosen as $(\sqrt{8n},4)$, and the rescaling of a variable will not influence the parameter $a$. We state the following concentration result, which could be found from \cite{boucheron2013concentration}.
\begin{theorem}
\label{thm:2.1}
For sub-exponential random variable X with parameter (a,b),
\begin{equation}
\label{sub}
P(X-EX\ge t)\le e^{max\{-\frac{t^2}{2a^2},-\frac{t}{2b}\}}
\end{equation}
\end{theorem}

Here we donote $X=x'^TA^TAx'$ and $Y=x^TA^TAx$, then $EX=EY=n$, for some $\varepsilon$
\begin{equation*}
\begin{aligned}
P(X-Y\ge -\varepsilon n^{\frac{3}{4}})&=P(Y-X\le\varepsilon n^{\frac{3}{4}})\\
 &=1-P(Y-X\ge \varepsilon n^{\frac{3}{4}})\\
 \eqref{sub}&\ge 1-e^{\max\{-\frac{\varepsilon n^{\frac{3}{4}}}{8},-\frac{\varepsilon^2 n^{\frac{1}{2}}}{16}\}}
\end{aligned}
\end{equation*}
So with high probability,
\begin{equation}
\label{f-f}
\begin{aligned}
f(\beta)-f(\beta')=&\frac{1}{n}(\sum\limits_{i=1}^n y_i\langle {x_i},{\beta}\rangle-\sum\limits_{i=1}^n y_i\langle {x_i},{\beta'}\rangle)+\frac{1}{2n}(Y-X)\\
\ge &\frac{\lambda\delta}{2}-\frac{\sqrt{Cj^2s\log(\frac{2p}{j^2s})}+4\gamma}{\sqrt{n}}-\frac{\varepsilon}{2\sqrt[4]{n}}\\
:=&\frac{\lambda\delta}{4}-\frac{a}{\sqrt{n}}-\frac{\varepsilon}{2\sqrt[4]{n}}.
\end{aligned}
\end{equation}
There, for convenience's sake, we set $\gamma=\frac{\lambda\delta\sqrt{n}}{8}$.\\

When $n>( \frac{\varepsilon+\sqrt{\varepsilon^2+4\lambda\delta a}}{\lambda\delta})^4$, the expression in \eqref{f-f} $>0$, so such $x' $ could not be a optimizer since $x$ is also a feasible point. Thus, our original  assumption that $\|\beta-\beta'\|_2\ge \delta$ cannot hold. This proves what we want.\\

Now we calculate the so-called ``high probability". We take $\varepsilon=\sqrt{\lambda\delta a}$ here. This statement holds when both parts in \eqref{plan2013} and \eqref{sub} hold. So P$\ge 1-4e^{-2\gamma^2}-e^{\max\{-\frac{\varepsilon n^{\frac{3}{4}}}{8},-\frac{\varepsilon^2 n^{\frac{1}{2}}}{16}\}}$. It is sufficient to set 
\[n\ge 4(\frac{\varepsilon+\sqrt{\lambda\delta a}}{\lambda\delta})^4=C\lambda^{-2}\delta^{-2} a^2\]
 and in this case, we have
\begin{equation*}
\begin{aligned}
P\ge &1-4e^{-2\gamma^2}-e^{\max\{-c\gamma^{\frac{3}{2}},-c'\gamma\}} \\
\ge &1-4e^{-2\gamma^2}-e^{-c\gamma}\\
\ge& 1-5e^{-c\gamma}\\
=&1-5e^{-c\delta \lambda\sqrt{n}}.
\end{aligned}
\end{equation*}

\subsection{Proof of Theorems \ref{thm:1.3}}
Second, we come to Theorem \ref{thm:1.3}.\\
In this section, we abbreviate $\hat{\beta_1}$ as $\hat \beta (\frac{2}{\lambda},1)$, and $\hat{\beta_k}$ as $\hat \beta (1,k)$.
First, notice that $\frac{1}{k}\hat{{\beta_k}}$ is feasible for problem
\[\min_\beta ||{y}-{X\beta}||_2\hspace{1cm}s.t.\ \|\beta\|_1\le \frac{2\sqrt{s}}{\lambda} ,\|\beta\|_2= 1.\]
and $k \hat{{\beta_1}}$ is feasible for 
\[\min_\beta ||{y}-{X\beta}||_2\hspace{1cm}s.t.\ \|\beta\|_1\le \sqrt{s},\|\beta\|_2= k.\]
So, from the definition of optimizer, we can conclude 
\begin{equation}
f(\hat{{\beta_k}})\ge f(k\hat{{\beta_1}})
\end{equation}
\begin{equation}
f(\hat{{\beta_1}})\ge f(\frac{1}{k}\hat{{\beta_k}})
\end{equation}

Denote $a=\frac{1}{n}\sum\limits_{i=1}^n y_i\langle {x_i},{\frac{\hat{{\beta_k}}}{k}}\rangle$ and $b=\frac{1}{2n}\sum\limits_{i=1}^n \langle {x_i},{\frac{\hat{{\beta_k}}}{k}}\rangle^2$ here. And, we similarly define $a'$, $b'$ for $\hat{{\beta_1}}$. Now $f(\hat{{\beta_k}})=ka-k^2b$, $f(k\hat{{\beta_1}})=ka'-k^2b'$.\\

Then the conditions shown above can be translated to
\[k a -k^2 b \ge k a' -k^2 b'\]
\[a'-b'\ge a-b \]
We deduce as following:
\begin{equation}
\label{ab}
\begin{aligned}
a-b&=a-kb +(k-1)b\\
 &\ge a'-kb'+(k-1)b\\
 &=a'-b'+(1-k)(b'-b)
\end{aligned}
\end{equation}
Which indicates
\[0\le (a'-b')-(a-b)\le (1-k)(b-b')\]
\[\Rightarrow \hspace{0.3cm}b\ge b'\]

Since $b-b'$ is sub-exponential, it cannot be too large, which indicates $\|f(\hat{\beta_1})-f(\frac{\hat{\beta_k}}{k})\|_2$ could not be too large. Thus, the same as Section \ref{secpro2}, $\hat{\beta_1}$ and $\frac{\hat{\beta_k}}{k}$ cannot be far away. We take this idea down to the ground. \\

According to \eqref{sub} and \eqref{f-f}, with $P\ge 1-e^{\max\{-\frac{\varepsilon n^{\frac{3}{4}}}{8},-\frac{\varepsilon^2 n^{\frac{1}{2}}}{16}\}}$, we have 
\begin{equation}
\label{b-b}
b-b' \le \frac{\varepsilon}{2\sqrt[4]{n}}
\end{equation}
When $\|\hat{\beta_1}-\frac{\hat{\beta_k}}{k}\|\ge \delta$, as shown in Section \ref{secpro2}, with P$\ge1-5e^{-c\delta \lambda\sqrt{n}}$
\[\|{\beta}-\hat{{\beta_1}}\|_2 \le \frac{\delta}{2}\]
\[\Rightarrow \|{\beta}-\frac{\hat{{\beta_k}}}{k}\|_2 \ge \frac{\delta}{2}\]
According to \eqref{f-f}
\begin{equation*}
f(\beta)-f(\frac{\hat{\beta_k}}{k})\ge \frac{\lambda\delta}{8}-\frac{a}{\sqrt{n}}-\frac{\varepsilon}{2\sqrt[4]{n}}
\end{equation*}
We choose 
\begin{equation}
\label{ep}
\varepsilon=\min\{\frac{\lambda}{2},\frac{1-\lambda}{2}\} \le \frac{\lambda\delta}{8}\sqrt[4]{n}-\frac{a}{\sqrt[4]{n}}
\end{equation}
combined with \eqref{b-b}
\begin{equation}
\frac{\lambda\delta}{8}-\frac{a}{\sqrt{n}}-\frac{\varepsilon}{2\sqrt[4]{n}}\ge\frac{\varepsilon}{2\sqrt[4]{n}}>(1-k)(b-b')
\end{equation}

Thus the assumption of $\|\hat{\beta_1}-\frac{\hat{\beta_k}}{k}\|\ge \delta$ could not hold with P$\ge1-5e^{-c\delta \lambda\sqrt{m}}$. On the inverse overally, we can conclude $\|\hat{\beta_1}-\frac{\hat{\beta_k}}{k}\|\le \delta$ with probability
\begin{equation*}
\begin{aligned}
P\ge & 1-e^{\max\{-\frac{\varepsilon n^{\frac{3}{4}}}{8},-\frac{\varepsilon^2 n^{\frac{1}{2}}}{16}\}}-5e^{-c\delta \lambda\sqrt{n}}\\
\ge & 1-6e^{-c\delta \lambda\sqrt{n}}
\end{aligned}
\end{equation*}

\subsection{Proof of Theorems \ref{thm:1.4}}
\label{secpro4}

Here in Section \ref{secpro4}, for every $\beta'$, we denote $\beta'=\mu {p}$, where $\mu=\|\beta'\|_2$, and $g(\mu)=f(\beta')$ for a fixed $p$.\\

First, we prove $\|\hat{{x_0}}\|_2 \le \lambda +\varepsilon$. The basic idea of our proof is that, since the optimal function is quadratic, it has a maximum point. When searching for the overall max point, the 2-norm of it will spontaneously come to a given point-$\lambda$. We would like to search for another point with larger function value if not so.\\
Notice $g'(\mu)=\frac{1}{n}(\sum\limits_{i=1}^n y_i\langle {x_i},{p}\rangle-\mu\sum\limits_{i=1}^n \langle {x_i},{p}\rangle^2)$. When $\mu >\lambda+\varepsilon$, noticing that $ng'(\mu)$ is sub-exponential, we conclude through concentration inequality:
\begin{equation*}
\begin{aligned}
P(g'(\mu)\le -c)&=P(ng'(\mu)\le -cn)\\
&=1-P(ng'(\mu)\ge -cn)\\
&=1-P\big(ng'(\mu)\ge E(g'(\mu)) -cn-E(\sum\limits_{i=1}^n y_i\langle {a_i},{p}\rangle-\mu\sum\limits_{i=1}^n \langle {a_i},{p}\rangle^2)\big)\\
&=1-P\big(ng'(\mu)\ge E(g'(\mu)) -cn-(\lambda-\mu)n\big)\\
&=1-P(ng'(\mu)-E(g'(\mu))\ge (\mu -\lambda -c)n)\\
&\ge 1-e^{\max\{-\frac{(\mu-\lambda-c)^2n}{8},-\frac{(\mu-\lambda-c)n}{8}\}}\\
&=1-e^{-\frac{(\mu-\lambda-c)^2n}{8}}\\
&\ge 1-e^{-\frac{\varepsilon^2n}{32}}
\end{aligned}
\end{equation*}
In the final row, we actually choose $c=\frac{\varepsilon}{2}$.\\

With this probability, for any $\beta'$ with $\|\beta'\|_2\ge \lambda+\varepsilon$, we can find another new feasible point through the line: $y=t\beta'$. Specifically, we can find $\bar{\beta}=\mu' {p}$, where $\mu'<\mu$, to satisfy $\bar{\beta}$ satisfies the restriction and $f(\bar{\beta})>f(\beta')$. The feasibility of $\bar{\beta}$ is quite obvious. So we have $\|\hat{{\beta_0}}\|_2 \le \lambda +\varepsilon$.\\

On the other hand, we prove $\|\hat{{\beta_0}}\|_2 \ge \lambda -\varepsilon$. This time, we go to the following proposition
\[with\ P\ge\ 1-e^{\frac{(1-\lambda-\varepsilon)^2\varepsilon^2}{8}n},\ f(\hat{{\beta_\lambda}})>f(\beta),\ where\ \beta\ is\ the\ solution\ to\]
\[\beta  = \arg\min_{\beta'} \sum_{i=1}^n(y_i - x_i^T\beta')^2, \ \ \ s.t.\   \|\beta'\|_1 \le s,\ \|\beta'\|_2\le\lambda-\varepsilon \]
\[\ then,\ with\ P\ge 1-e^{\frac{(1-\lambda-\varepsilon)^2\varepsilon^2}{8}n},\ the\ solution\ of\ non-scale\ cannot\ lie\ into\ the\ disk\ of\ \lambda-\varepsilon\]
\vspace{0.2cm}

Assume $\|\beta\|_2=k,\ then\ k\le\lambda-\varepsilon$ and according to \eqref{sub}.
\begin{equation*}
\begin{aligned}
P(f(\hat{{\beta_\lambda}})-f(\beta)>c)&=P(f(\hat{{\beta_\lambda}})-f(\beta)-E(f(\hat{{\beta_\lambda}})-f(\beta))>c-E(f(\hat{{\beta_\lambda}})-f(\beta)))\\
&=1-P(f(\beta)-f(\hat{{\beta_\lambda}})-E(f(\beta)-f(\hat{{\beta_\lambda}}))\ge E(f(\hat{{\beta_\lambda}})-f(\beta))-c)\\
&=1-P(f(\beta)-f(\hat{{\beta_\lambda}})-E(f(\beta)-f(\hat{{\beta_\lambda}}))\ge Ef(\hat{{\beta_\lambda}})-Ef(\beta)-c)\\
&=1-P(f(\beta)-f(\hat{{\beta_\lambda}})-E(f(\beta)-f(\hat{{\beta_\lambda}}))\ge \lambda-k-\frac{1}{2}(\lambda^2-k^2)-c)\\
&=1-P(nf(\beta)-nf(\hat{{\beta_\lambda}})-E(nf(\beta)-nf(\hat{{\beta_\lambda}}))\ge (\lambda-k)(1-\frac{\lambda+k}{2})n-cn)\\
&\ge 1-e^{\max\{-\frac{((\lambda-k)(1-\frac{\lambda+k}{2})-c)^2}{8}n,-\frac{(\lambda-k)(1-\frac{\lambda+k}{2})-c}{8}n\}}\\
&= 1-e^{-\frac{((\lambda-k)(1-\frac{\lambda+k}{2})-c)^2}{8}n}\\
&\ge 1-e^{-\frac{(\varepsilon (2-2\lambda-\varepsilon)-2c)^2}{32}n}
\end{aligned}
\end{equation*}
Thus, we set $c=\frac{\varepsilon^2}{4}$ and conclude
\[P(f(\hat{{\beta_\lambda}})-f(\beta)>\frac{\varepsilon^2}{4})\ge 1-e^{-\frac{(1-\lambda-\varepsilon)^2\varepsilon^2}{8}n}\]
So, with this probability, the optimizer of non-scale cannot lie inside the $\lambda-\varepsilon$ disc. Finally, we combine these two parts to calculate the overall probability
\[P\ge 1-e^{-\frac{\varepsilon^2n}{32}}-e^{-\frac{(1-\lambda-\varepsilon)^2\varepsilon^2}{8}n}.\]

Till now, all proof is finished.$\qed$

\section{Simulation}
\label{sec:simu}
In this section, we use simulation study to illustrate our theoretical findings. We first demonstrate that when data are from a logistic model which is a special single index model by choosing $F(t) = \frac{e^t}{1+e^t}$ in model \eqref{eqn:single-index} .  In this special model, the classifier is defined as
 \[\ell(x) = sign(x^T\hat\beta),\] where $\hat\beta$ could be estimated via vanilla Lasso -- $\ell_1$ constrained least squares. For a new data point with observed features $x$, we calculate $\ell(x)$ and claim it as labelled +1 if $\ell(x) = 1$ or else, we claim it as labelled -1.  From this classifier, we see that the scale $\hat\beta$ does not make any contribution to the classifier. In fact, we will show later that $\beta^*$ could not be estimated correctly, but $\frac{\beta^*}{\|\beta^*\|_2}$ could.

Let $\beta^*$ be a fixed s-sparse vector in $\mathbb R^p$, and ${X}_{n\times p}$ is an $n\times p$ matrix with entries from i.i.d Gaussian distribution where $n$ is the number of samples. Note that in this section, we set $\beta^*$ as any s-sparse vector, so the 2-norm of it is probably not 1. We generate a response vector $y$ from the following model:
\[y_i\in \{-1,1\}\]
\[P(y_i=1|x_i)=\frac{1}{1+\exp\{x_i^T\beta\}}\]
\[P(y_i=-1|x_i)=1-P(y_i=1|x_i)\]
Here, we set $p=1200,\ s=10$, and $n$ from 200 to 3000. Matrices ${X}$ are generated for each $n$. Results are showed below, in which $x$-axis stands for $200n$ and $y$-axis for error $||\frac{\beta^*}{\|\beta^*\|_2} -\frac{\hat\beta}{||\hat{\beta}||_2}||_2$:

\begin{figure}[h!]
\begin{center}
\includegraphics[scale=0.5]{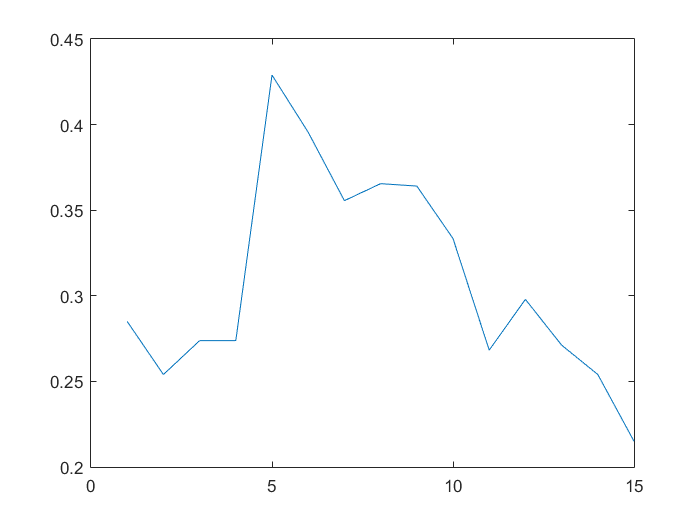}
\caption{Relationship of $n$ and error of $\beta$. \label{fig:error}}
\end{center}
\end{figure}
As we can see, the error declines fairly sharply as sample size $n$ grows.  We also show the differences between $\hat\beta$ and $\beta^*$. From Figure \ref{fig:error} we see that even when $n$ is very big, the difference does not vanish, which suggests that $\hat\beta$ is not a good estimator for $\beta^*$, while Figure \ref{fig:error} suggests that $\hat\beta/\|\hat\beta\|_2$ is a good estimator of $\beta^*/\|\beta^*\|_2$.

\section{Conclusions}
\label{sec:conclusion}
In this paper, we studied how Lasso performs for a sparse classification problem when data from a single-index model.  Vanilla Lasso is very simple and is computationally very fast. For high-dimensional data, it is much more convenient than using semi-parametric method directly. We show that it performs quite well for high-dimensional binary classification problems. If one classifier like logistic regression or probit regression model does not depend on the scale of $\beta$, vanilla Lasso could give direct classifiers. If one does want to use more complex models, vanilla lasso could be used as a first step to remove huge amount of redundant features.

\bibliographystyle{chicago}
\bibliography{ref}
\end{document}